\documentclass[10pt]{article}
\usepackage{amsmath,amsthm,amsfonts,amssymb,latexsym,graphicx}
\usepackage[noadjust]{cite}
\usepackage{algorithm}
\usepackage[noend]{algpseudocode}
\algrenewcommand\algorithmicthen{\relax}
\algrenewcommand\algorithmicdo{\relax}
\usepackage[pdfpagemode=UseNone,pdfstartview=FitH]{hyperref}

\usepackage{booktabs}
\usepackage{cellspace}
\allowdisplaybreaks[4]

  \theoremstyle{definition}

\DeclareMathOperator{\UB}{UB}
\DeclareMathOperator{\LB}{LB}

\DeclareMathOperator{\Markov}{Markov}
\DeclareMathOperator{\Ber}{Ber}

\title{Conformal testing: binary case with Markov alternatives}
\author{Vladimir Vovk, Ilia Nouretdinov, and Alex Gammerman}

\begin{document}
\maketitle

\begin{abstract}
  We continue study of conformal testing in binary model situations.
  In this note we consider Markov alternatives to the null hypothesis of exchangeability.
  We propose two new classes of conformal test martingales;
  one class is statistically efficient in our experiments,
  and the other class partially sacrifices statistical efficiency
  to gain computational efficiency.

  The version of this paper at \url{http://alrw.net} (Working Paper 36)
  is updated most often.
\end{abstract}

\section{Introduction}

This note treats a problem similar to the one considered in \cite{Ramdas/etal:arXiv2102}:
we would like to test online the null hypothesis of exchangeability of binary observations
under Markov alternatives.

The simplest way of online hypothesis testing is to use \emph{test martingales},
which are defined as nonnegative processes with initial value 1
that are martingales under the null hypothesis;
see, e.g., \cite{Shafer/etal:2011}.
Such processes, for the null hypothesis of exchangeability,
can be constructed using the method of conformal prediction \cite{Vovk/etal:2005book},
and we will refer to them as conformal test martingales.
A previous paper \cite{Vovk:2021COPA} constructs custom-made conformal test martingales
for different alternative hypotheses, those of a changepoint.

The method of \cite{Ramdas/etal:arXiv2102}, which is specifically devoted to Markov alternatives,
is more general:
instead of a test martingale the authors construct a ``safe e-pro\-cess'' (to be defined in the next section).
Safe e-processes are closely related to test martingales and admit a similar interpretation
as the capital of a gambler trying to discredit the null hypothesis.
Our methods give similar results to the methods of \cite{Ramdas/etal:arXiv2102} in the model situations that we consider
(following \cite{Ramdas/etal:arXiv2102}).
The advantage of our methods is that they extend easily to the usual setting of machine learning,
where the observations are pairs $(x,y)$ consisting of a potentially complex object $x$ and its label $y$.

In this note we only design conformal test martingales for a simple alternative hypothesis
(a specific probability measure).
This is different from \cite{Ramdas/etal:arXiv2102},
who are interested in testing against the composite alternative Markov hypothesis.
As in \cite{Ramdas/etal:arXiv2102},
we could mix our conformal test martingales over the possible alternative hypotheses,
but we leave this step for future research.

\section{Model situations}

\begin{figure}[b]
  \begin{center}
    \includegraphics[width=0.48\textwidth]{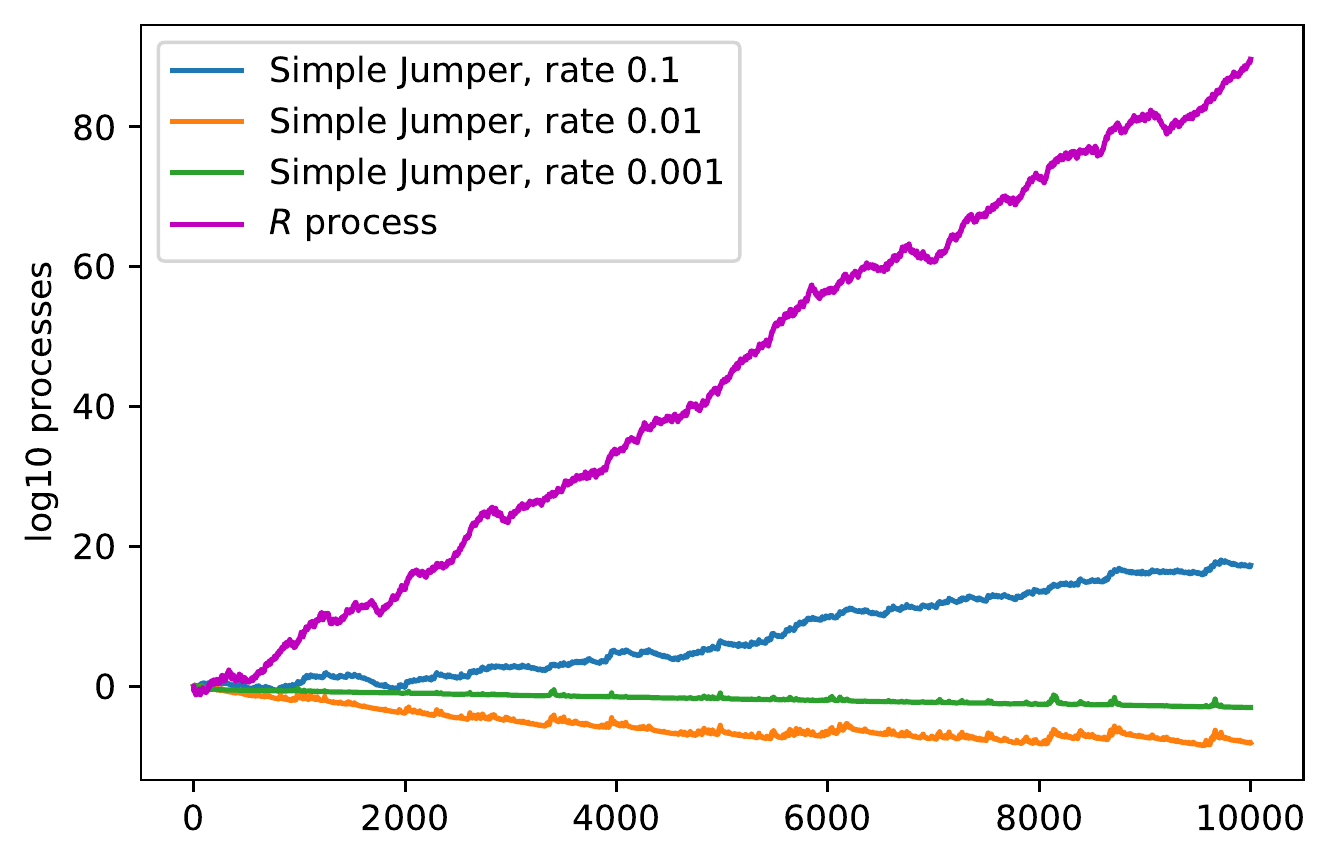}
    \includegraphics[width=0.48\textwidth]{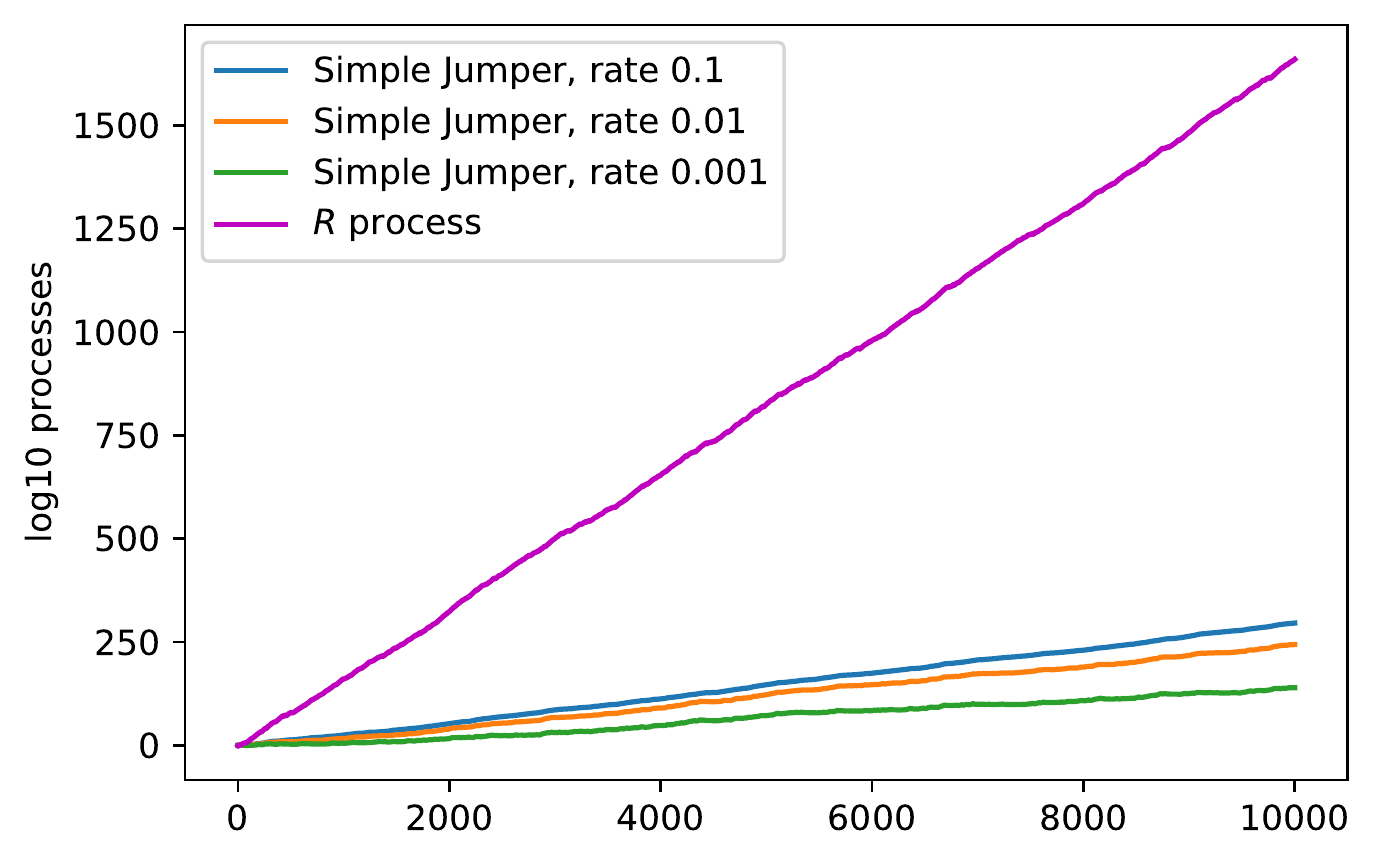}
  \end{center}
  \caption{The process $R$ of \cite{Ramdas/etal:arXiv2102} and the Simple Jumper in the large scenario.
    Left panel: the hard case.
    Right panel: the easy case.}
  \label{fig:L_existing}
\end{figure}

This section introduces the model situations considered in this paper,
following \cite[Section 4.2]{Ramdas/etal:arXiv2102}.
Our data consist of binary observations generated from a Markov model.
We will use the notation $\Markov(\pi_{1|0},\pi_{1|1})$
for the probability distribution of a Markov chain with the transition probabilities
$\pi_{1|0}$ for transitions $0\to1$ and $\pi_{1|1}$ for transitions $1\to1$;
the probability that the first observation is 1 will always be assumed $0.5$.
In the \emph{hard case}, the model is $\Markov(0.4,0.6)$,
and in the \emph{easy case}, the model is $\Markov(0.1,0.9)$.
The number of observations is $N:=10^4$ (as in \cite{Ramdas/etal:arXiv2102}) or $N:=10^3$ or $N:=10^2$;
we will refer to these scenarios as \emph{large}, \emph{medium}, and \emph{small}, respectively.

In all our experiments we use 2021 as the seed for the NumPy pseudorandom number generator.
(This, however, does not make the trajectories in our plots comparable between different scenarios.)
The dependence on the seed will be explored in boxplots reported in Section~\ref{sec:simplified};
the seed affects not only the data but also the values of conformal martingales,
which are randomized processes, given the data.

Let $B_{\pi}$ be the Bernoulli distribution on $\{0,1\}$ with parameter $\pi\in[0,1]$:
$B_{\pi}(\{1\})=\pi$.
Set $\Ber(\pi):=B_{\pi}^{\infty}$.
Our null hypothesis is the \emph{IID model},
under which the observations are generated from $\Ber(\pi)$ with unknown parameter $\pi$.

Ramdas et al.\ construct a \emph{safe e-process} $R=R_n$:
namely, under any $\Ber(\pi)$,
$R$ is dominated by a test martingale $M^{(\pi)}_n$ w.r.\ to $\Ber(\pi)$,
in the sense that $R_n\le M_n^{(\pi)}$ for all $n$ and $\pi$.
The trajectories of their process for the two cases, hard and easy, are shown in Figure~\ref{fig:L_existing}
(they coincides with those in Figure~4 in \cite{Ramdas/etal:arXiv2102}
apart from using base 10 logarithms and a different randomly generated dataset).
The figure also shows trajectories of the Simple Jumper martingale
(see, e.g., \cite{Vovk:2021})
for various values of the jumping rate; it performs poorly in this context.

\section{Two benchmarks}

In this section we will discuss possible benchmarks that we can use for evaluating the quality of our conformal test martingales.
The \emph{upper benchmark} is
\begin{equation*}
  \UB_n
  :=
  \frac{\Markov(\pi_{1|0},\pi_{1|1})([z_1,\dots,z_n])}{\Ber(0.5)([z_1,\dots,z_n])},
\end{equation*}
where $[z_1,\dots,z_n]$ is the set of all infinite sequences of binary observations
starting from $z_1,\dots,z_n$,
and $z_1,z_2,\dots$ are the actual observations.
The \emph{lower benchmark} is
\[
  \LB_n
  :=
  \frac{\Markov(\pi_{1|0},\pi_{1|1})([z_1,\dots,z_n])}{\Ber(\hat\pi)([z_1,\dots,z_n])},
\]
where $\hat\pi:=k/n$ (the maximum likelihood estimate) and $k=k(n)$ is the number of 1s among $z_1,\dots,z_n$.
By definition, $\UB_0=\LB_0:=1$.

\begin{figure}
  \begin{center}
    \includegraphics[width=0.48\textwidth]{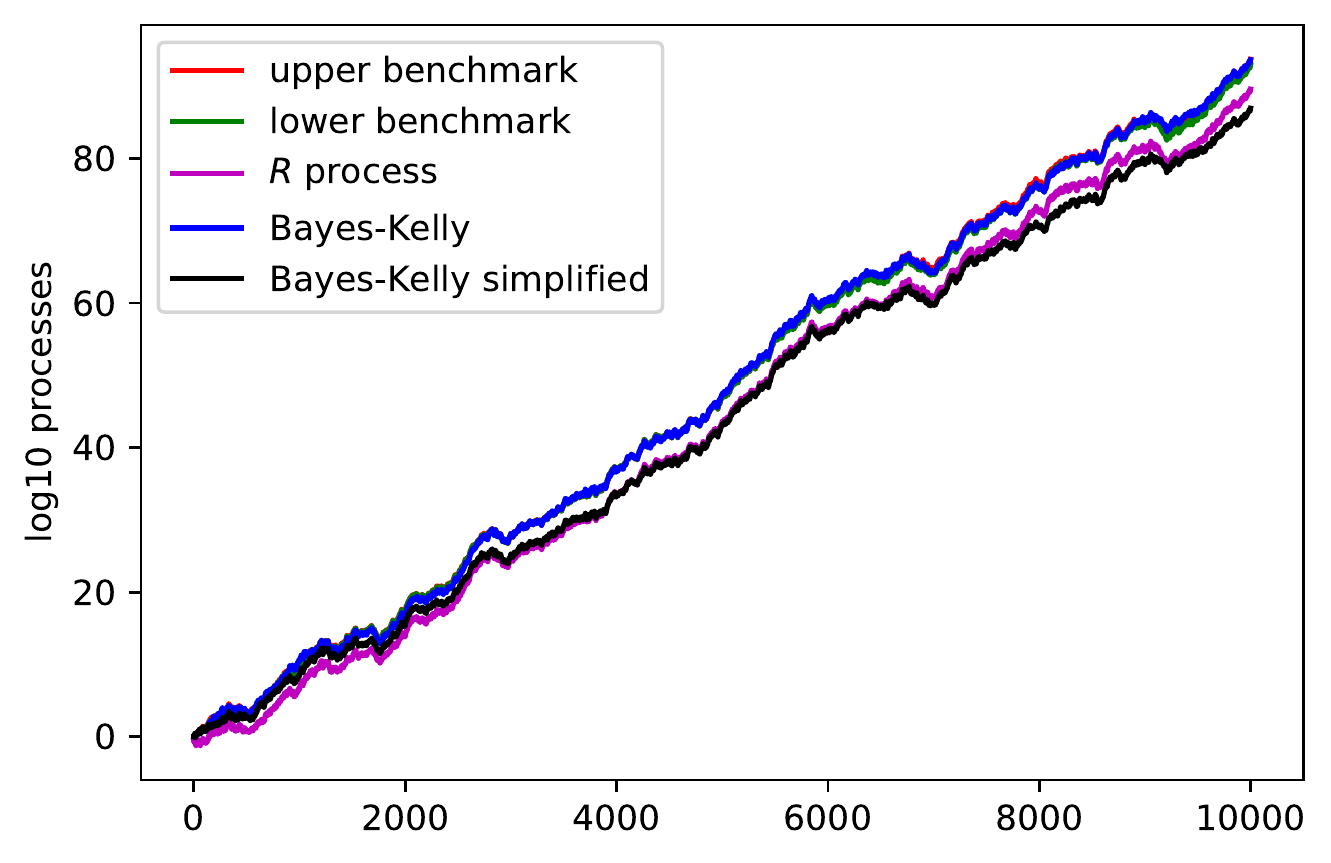}
    \includegraphics[width=0.48\textwidth]{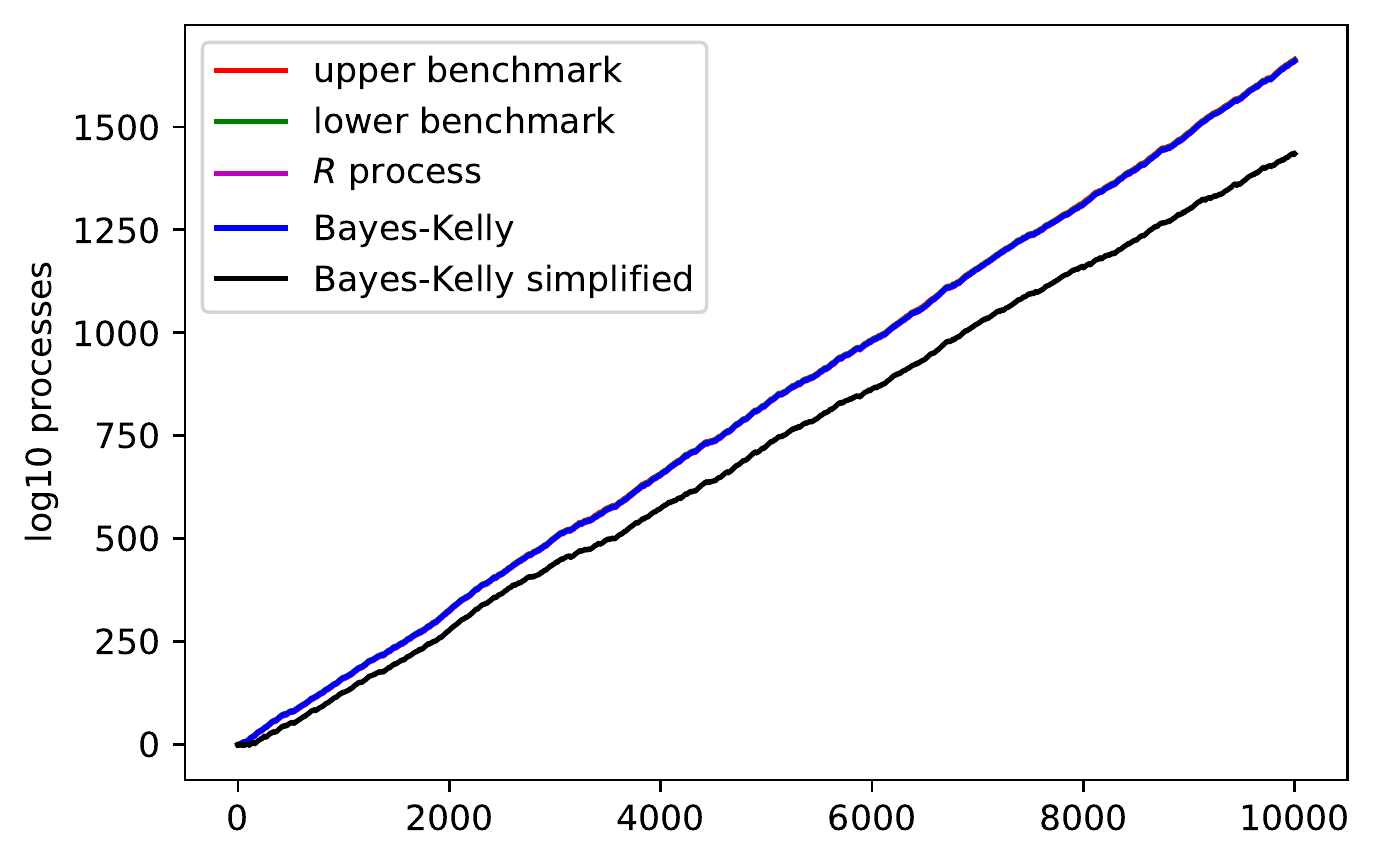}
  \end{center}
  \caption{The two benchmarks,
    $R$ process, Bayes--Kelly conformal test martingale, and its simplified version
    in the large scenario.
    Left panel: hard case.
    Right panel: easy case.}
  \label{fig:L_processes}
\end{figure}

\begin{figure}
  \begin{center}
    \includegraphics[width=0.48\textwidth]{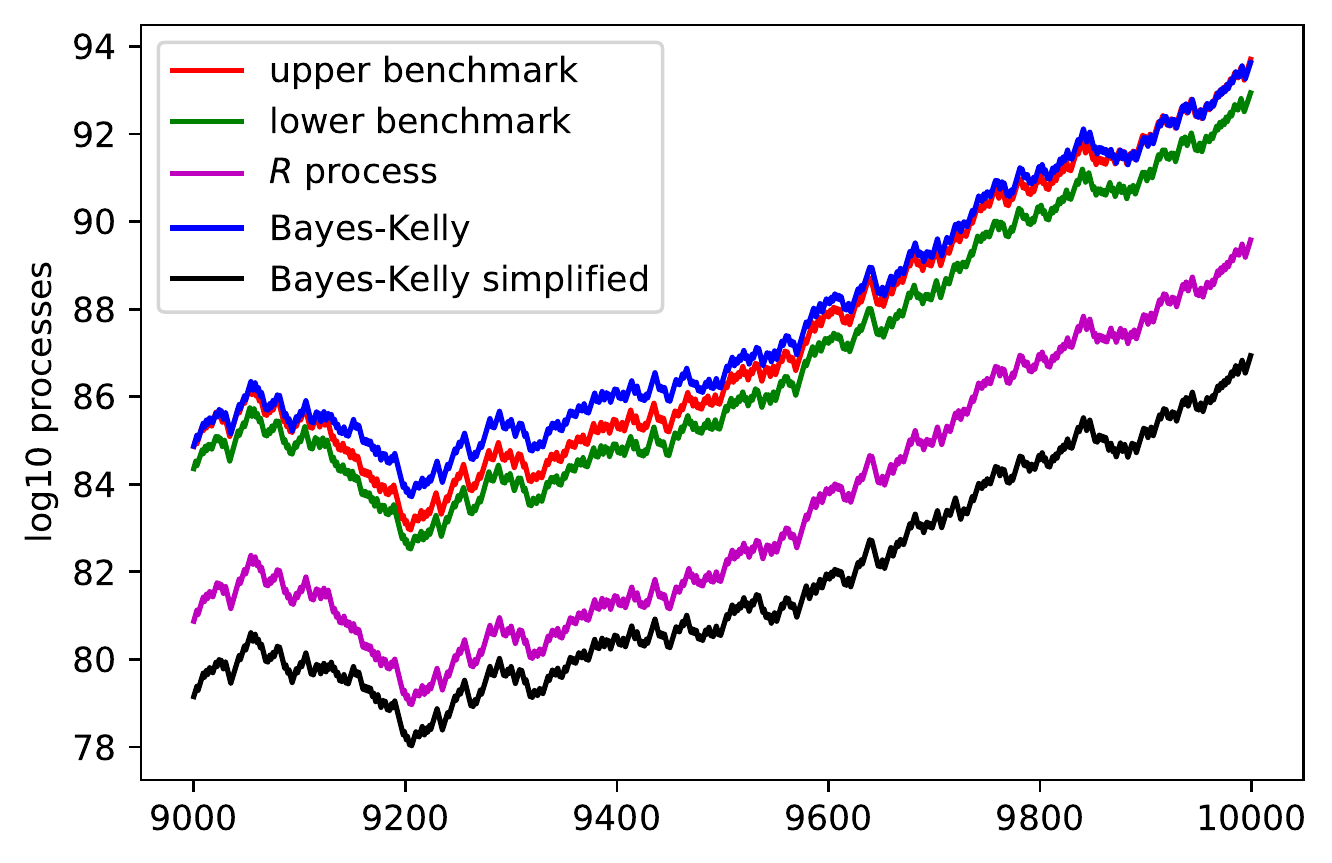}
    \includegraphics[width=0.48\textwidth]{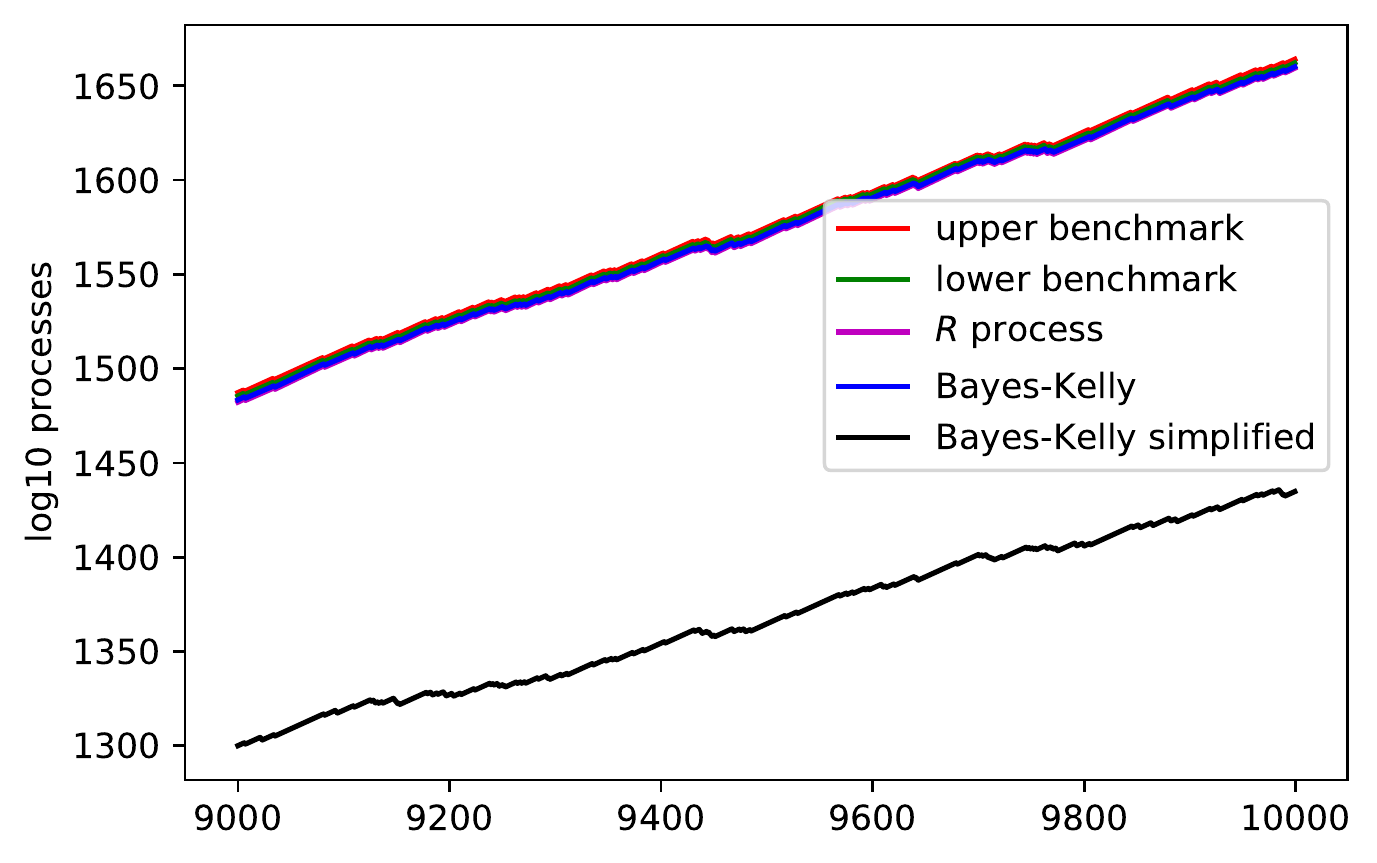}
  \end{center}
  \caption{The analogue of Figure~\ref{fig:L_processes} for the last 1000 observations.}
  \label{fig:Lend_processes}
\end{figure}

The trajectories of the upper and lower benchmarks are shown in Figure~\ref{fig:L_processes} in red and green;
the figure also shows the trajectory the $R$ process discussed in the previous section,
and the other two trajectories should be ignored for now.
The two benchmarks coincide or almost coincide.
Figure~\ref{fig:Lend_processes} should the same trajectories ``under the lens'',
over the last 1000 observations.

\section{Bayesian conformal testing}

In this section we will use a Bayesian method that is statistically efficient in our experiments
but whose computational efficiency will be greatly improved in the next section.
The p-values $p_1,p_2,\dots$ are generated as described in \cite{Vovk:2021COPA};
in particular, we are using the identity nonconformity measure
(the nonconformity score of an observation $z$ is $z$).
Under the alternative hypothesis, the p-values are generated by a completely specified stochastic mechanism.
According to \cite[Theorem 2]{Fedorova/etal:2012ICML},
the optimal (in the Kelly-type sense of that paper) betting functions $f_n$
are given by the density of the predictive distribution of $p_n$ conditional on knowing $p_1,\dots,p_{n-1}$.
Let us find these predictive distributions.
We will use the notation $U[a,b]$, where $a<b$, for the uniform probability distribution
on the interval $[a,b]$ (so that its density is $1/(b-a)$).

We are in a typical situation of Bayesian statistics.
The Bayesian parameter is the binary sequence $(z_1,z_2,\dots)\in\{0,1\}^{\infty}$ of observations,
and the prior distribution on the parameter is $\Markov(\pi_{1|0},\pi_{1|1})$.
The Bayesian observations are the conformal p-values $p_1,p_2,\dots$.
Given the parameter, the distribution of $p_n$ is
\[
  p_n
  \sim
  \begin{cases}
    U[0,k/n] & \text{if $z_n=1$}\\
    U[k/n,1] & \text{if $z_n=0$},
  \end{cases}
\]
where $k:=z_1+\dots+z_n$ is the number of 1s among the first $n$ observations.

Let $w^n_{k,j}$, where $n=1,2,\dots$, $k=0,\dots,n$, and $j\in\{0,1\}$,
be the total posterior probability of the parameter values $z_1,z_2,\dots$
for which $z_1+\dots+z_n=k$ and $z_n=j$;
we will use them as the weights when computing the predictive distributions for the p-values.
We can compute the weights $w^n_{k,j}$ recursively in $n$ as follows.
We start from
\[
  w^1_{0,0}
  :=
  w^1_{1,1}
  :=
  0.5,
  \quad
  w^1_{0,1}
  :=
  w^1_{1,0}
  :=
  0.
\]
At each step $n\ge2$,
first we compute the unnormalized weights
\begin{align*}
  \tilde w^n_{k,0}
  &:=
  \left(
    w^{n-1}_{k,0}
    \pi_{0|0}
    +
    w^{n-1}_{k,1}
    \pi_{1|0}
  \right)
  l^{n-1}_{k}(0,p_n),\\
  \tilde w^n_{k,1}
  &:=
  \left(
    w^{n-1}_{k-1,0}
    \pi_{1|0}
    +
    w^{n-1}_{k-1,1}
    \pi_{1|1}
  \right)
  l^{n-1}_{k-1}(1,p_n),
\end{align*}
where $l$ is the likelihood defined by
\begin{align*}
  l^n_k(1,p)
  &:=
  \begin{cases}
    \frac{n+1}{k+1} & \text{if $p\le\frac{k+1}{n+1}$}\\
    0 & \text{otherwise},
  \end{cases}\\
  l^n_k(0,p)
  &:=
  \begin{cases}
    \frac{n+1}{n-k+1} & \text{if $p\ge\frac{k}{n+1}$}\\
    0 & \text{otherwise},
  \end{cases}
\end{align*}
and then we normalize them:
\[
  w^n_k
  :=
  \tilde w^n_k / \sum_{k=0}^m\sum_{j=0}^1 \tilde w^n_{k,j}.
\]

Given the posterior weights for the previous step,
we can find the predictive distribution for $p_n$ as
\[
  p_n
  \sim
  \sum_{k=0}^{n-1}
  \sum_{j=0}^1
  w^{n-1}_{k,j}
  \left(
    \pi_{1|j} U\left[0,\frac{k+1}{n}\right]
    +
    \pi_{0|j}
    U\left[\frac{k}{n},1\right]
  \right),
\]
where we use the shorthand $\pi_{0|j} := 1-\pi_{1|j}$.
Therefore, the betting functions for the resulting \emph{Bayes--Kelly conformal test martingale} are
\begin{equation}\label{eq:BK}
  f_n(p)
  =
  \sum_{k=0}^{n-1}
  \sum_{j=0}^1
  w^{n-1}_{k,j}
  \left(
    \frac{n}{k+1}
    \pi_{1|j} 1_{p\le\frac{k+1}{n}}
    +
    \frac{n}{k}
    \pi_{0|j} 1_{p\ge\frac{k}{n}}
  \right).
\end{equation}

\begin{figure}
  \begin{center}
    \includegraphics[width=0.48\textwidth]{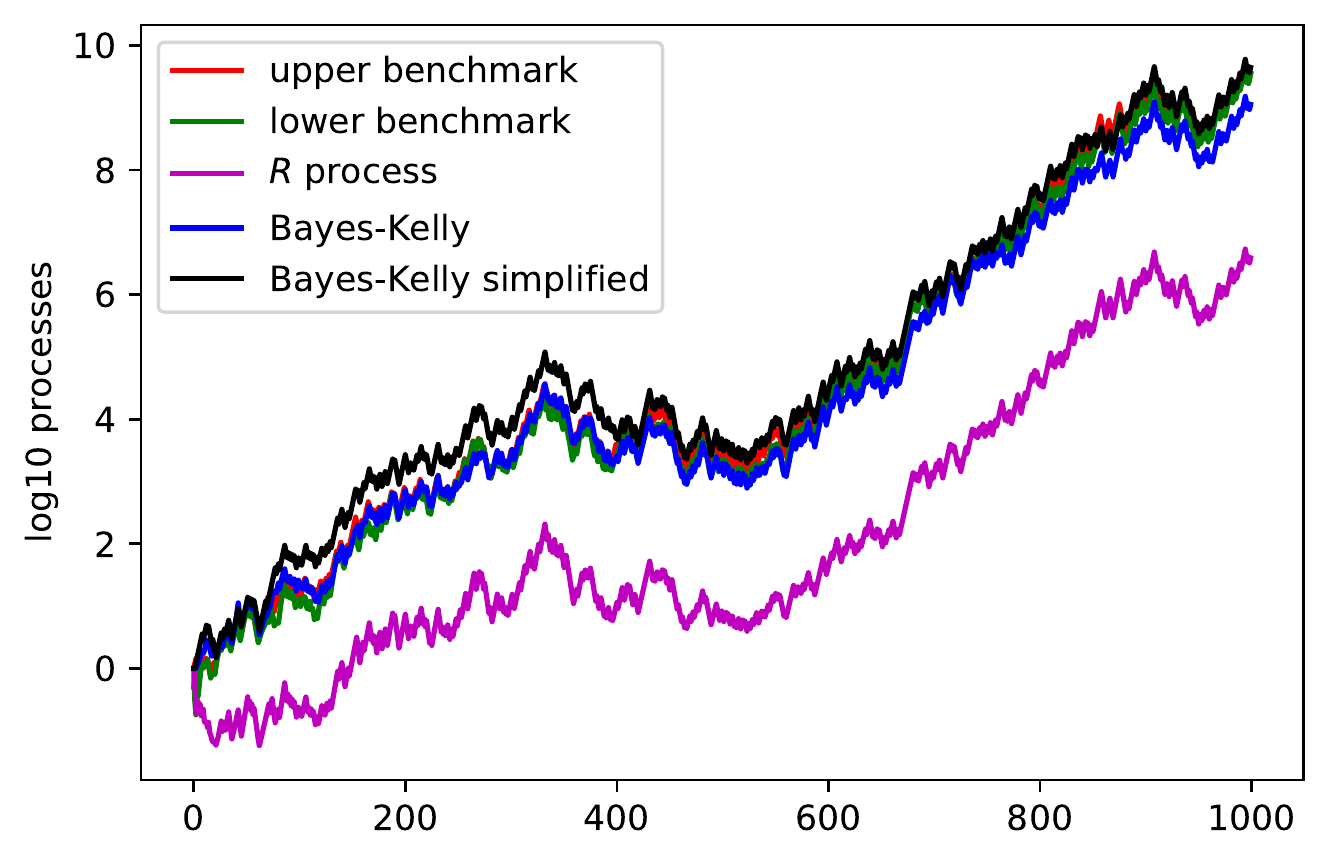}
    \includegraphics[width=0.48\textwidth]{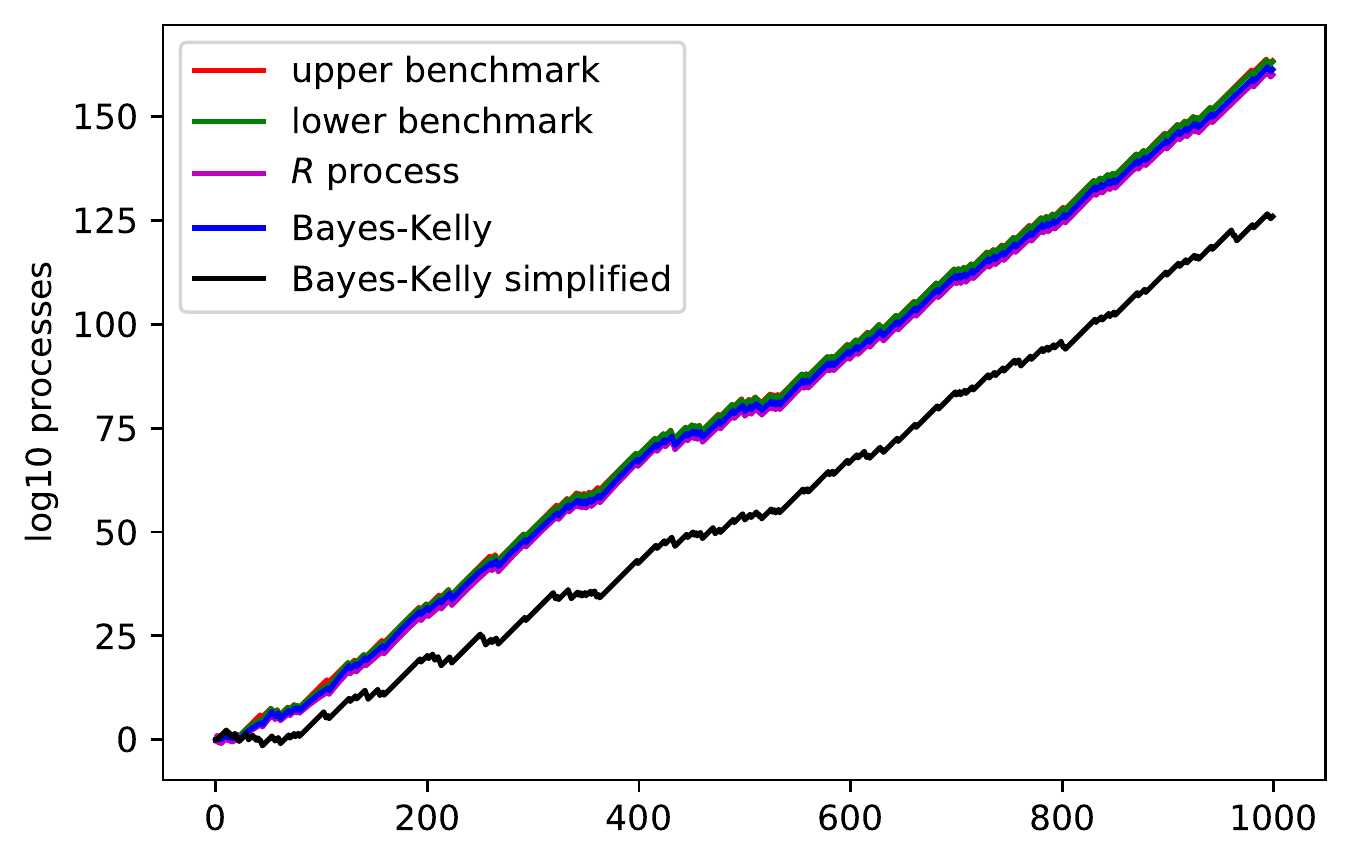}
  \end{center}
  \caption{The Bayes--Kelly and Bayes--Kelly simplified conformal test martingales, the $R$-process, and the two benchmarks
    in the middle scenario.
    Left panel: hard case.
    Right panel: easy case.}
  \label{fig:M_processes}
\end{figure}

For experimental results, see Figure~\ref{fig:M_processes},
in addition to Figure~\ref{fig:L_processes}.
The Bayes--Kelly conformal test martingale appears to be very close to the two benchmarks.
Its simplified version is described in the next section.
The relatively poor performance of the $R$-process in the left panel of Figure~\ref{fig:M_processes}
should not be interpreted as it being inferior to the Bayes--Kelly conformal test martingale:
remember that $R$ works against all Markov alternatives,
whereas the other processes in Figures~\ref{fig:L_processes}--\ref{fig:S_boxplots}
are adapted to the specific alternative hypothesis
($\Markov(0.4,0.6)$ in the hard case and $\Markov(0.1,0.9)$ in the easy case).

\section{Simplified Bayesian conformal testing}
\label{sec:simplified}

\begin{figure}
  \begin{center}
    \includegraphics[width=0.48\textwidth]{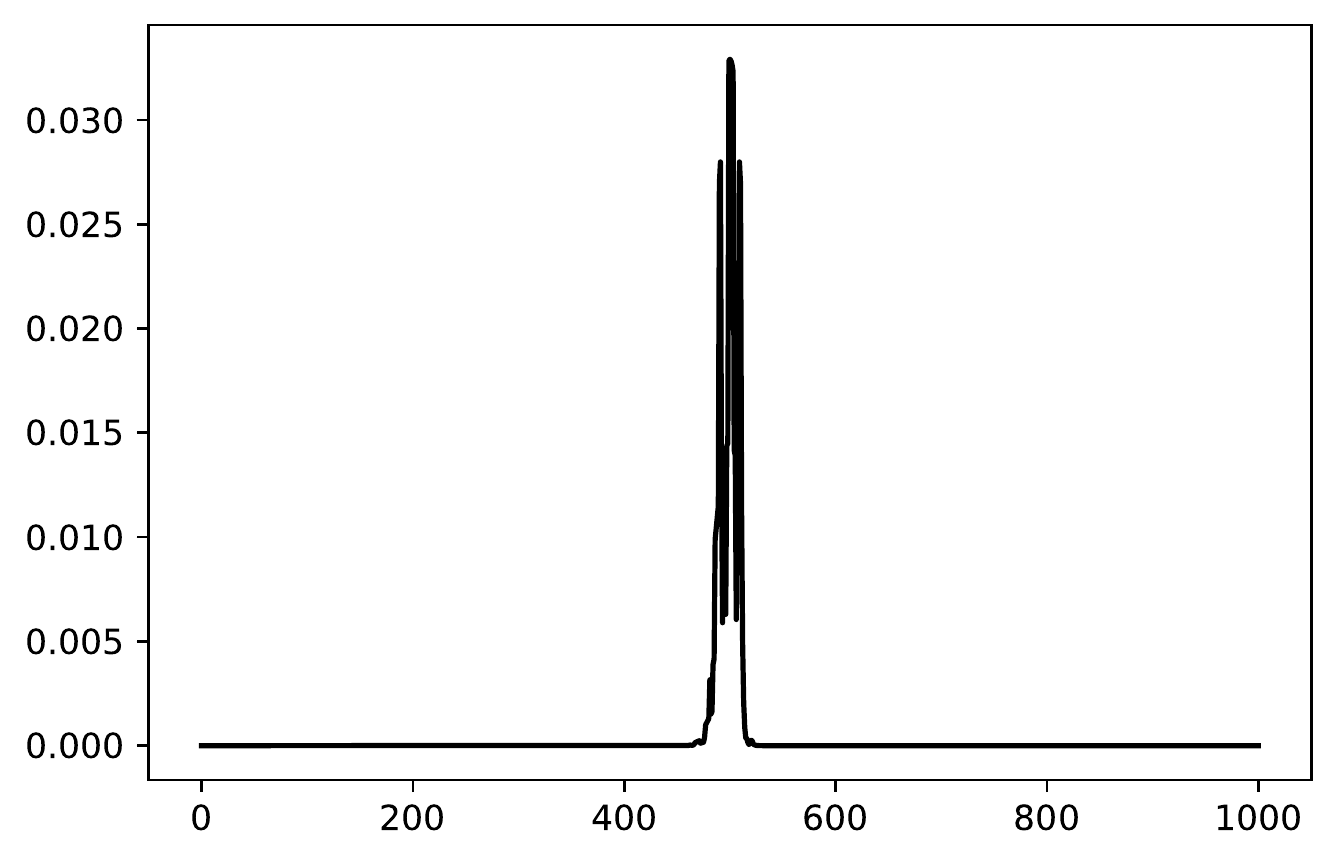}
    \includegraphics[width=0.48\textwidth]{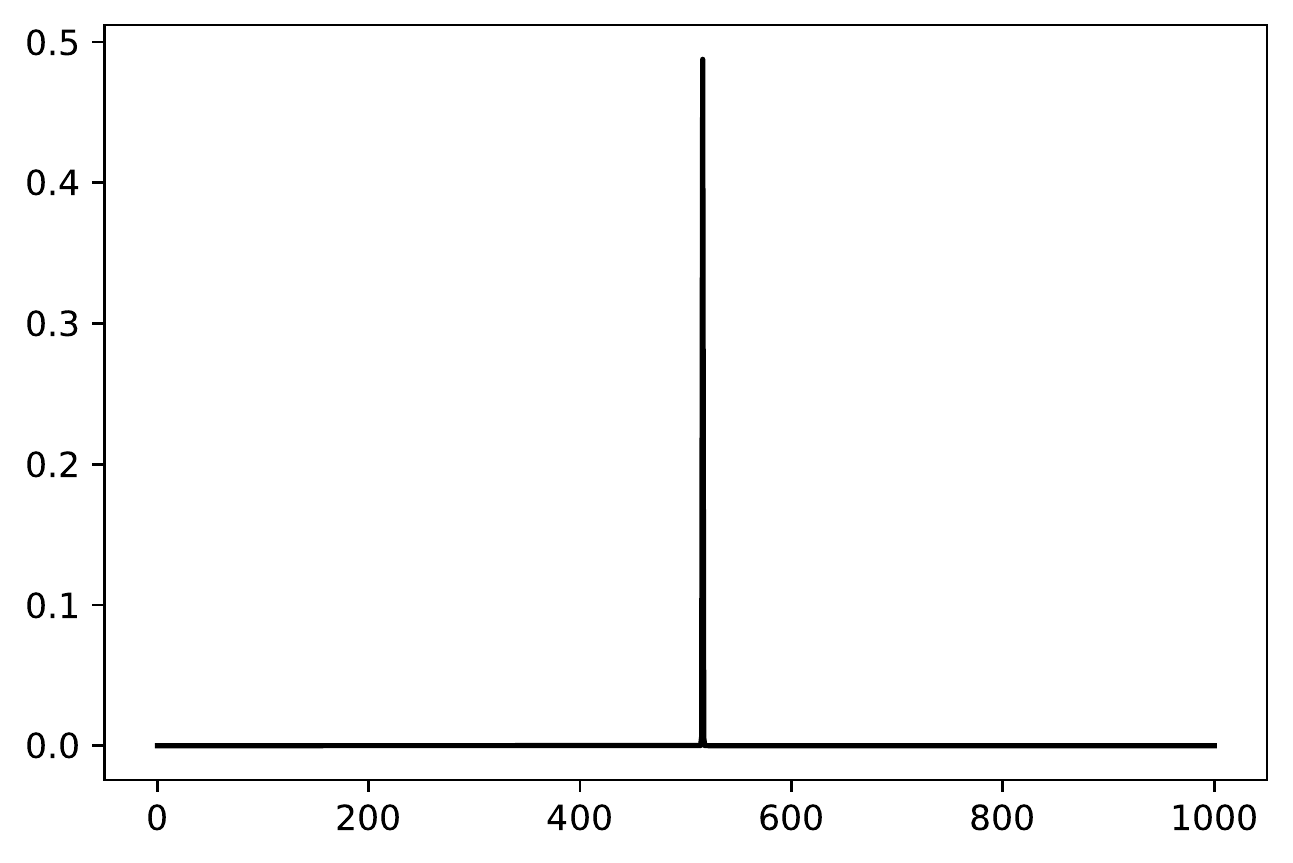}
  \end{center}
  \caption{The weights $w^{1000}_k$, $k=0,\dots,1000$,
    at the last step for the Bayes--Kelly conformal test martingale in the medium scenario
    (the hard case on the left and easy on the right).}
  \label{fig:M_weights}
\end{figure}

In this section we consider a radical simplification of the Bayes--Kelly conformal test martingale \eqref{eq:BK}.
We still assume that the Markov chain is symmetric, as in our model situations.
If we assume that the weights $w^n_{k,j}$, $k=0,\dots,n$,
are concentrated at
\[
  k \approx k+1 \approx n/2,
\]
\eqref{eq:BK} will simplify to
\begin{equation}\label{eq:BK-simplified}
  f_n(p)
  =
  2
  \pi_{1|j} 1_{p\le0.5}
  +
  2
  \pi_{0|j} 1_{p>0.5}.
\end{equation}
Figure~\ref{fig:M_weights} shows the weights (averaged over $j\in\{0,1\}$)
for the last step of the Bayes--Kelly conformal test martingale
in the medium scenario ($10^3$ observations).
They are indeed concentrated around values of $k$ not so different from $0.5 N = 500$.

\begin{figure}
  \begin{center}
    \includegraphics[width=0.48\textwidth]{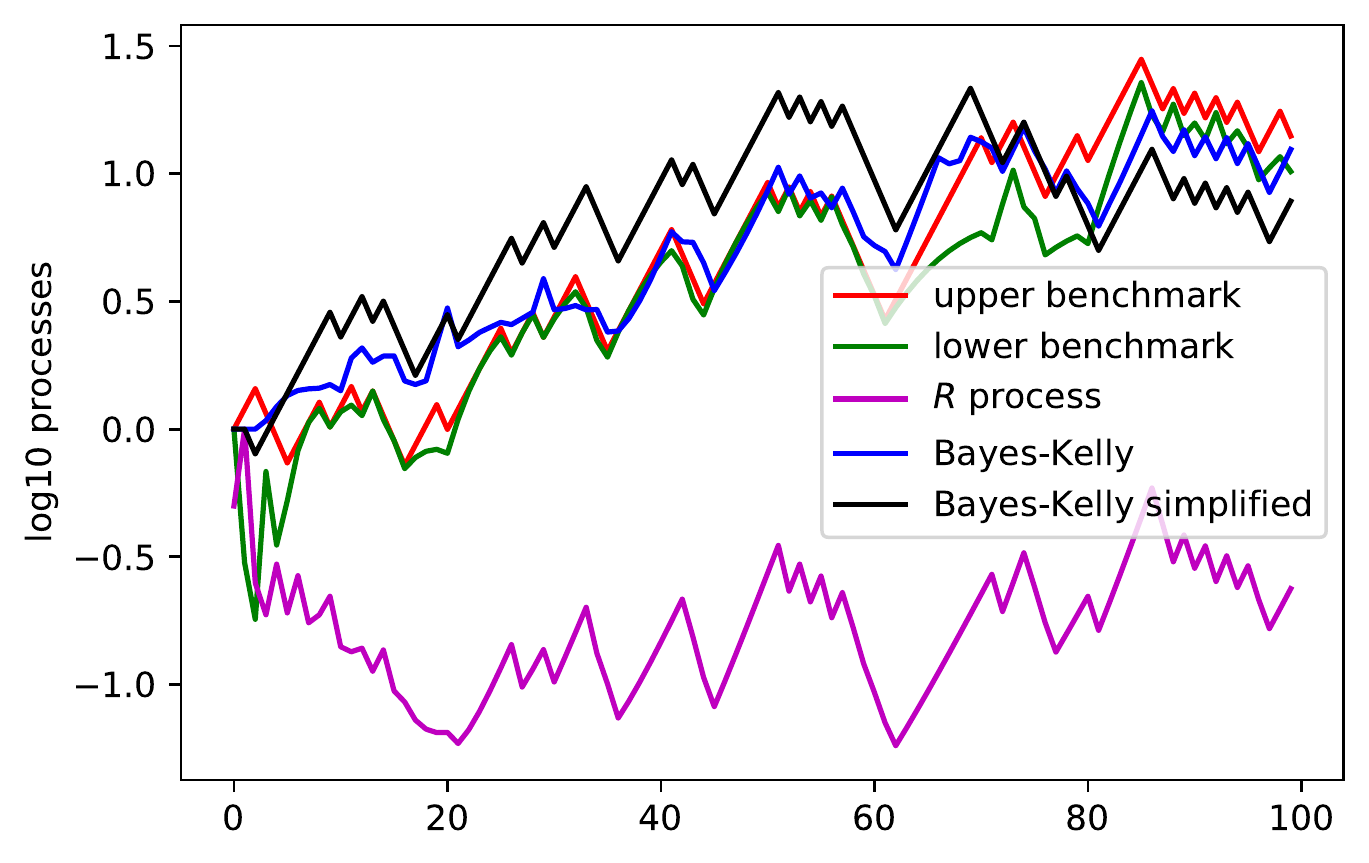}
    \includegraphics[width=0.48\textwidth]{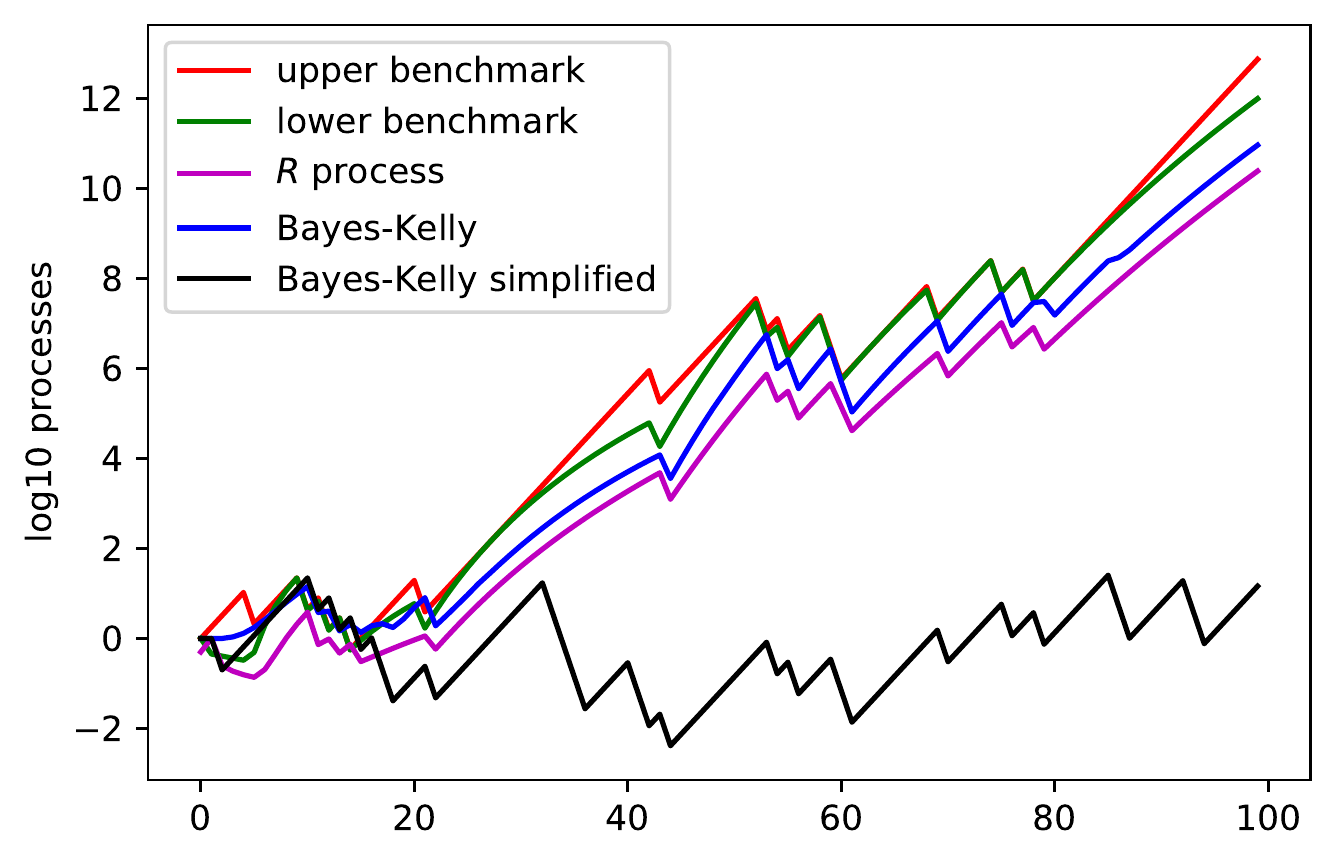}
  \end{center}
  \caption{The analogue of Figures~\ref{fig:L_processes} and~\ref{fig:M_processes}
    for the small scenario.}
  \label{fig:S_processes}
\end{figure}

As a second step, we make \eqref{eq:BK-simplified} straightforward to compute by setting
\[
  j
  :=
  \begin{cases}
    1 & \text{if $p_{n-1}\le0.5$}\\
    0 & \text{if not}.
  \end{cases}
\]
(If $k(n-1):=z_1+\dots+z_{n-1}\approx(n-1)/2$, then $j=z_{n-1}$ with high probability.)
The performance of the simplified version is shown
in Figures~\ref{fig:L_processes}--\ref{fig:M_processes} and \ref{fig:S_processes}.
It is usually worse than that of the Bayes--Kelly conformal test martingale and the two benchmarks,
but is comparable on the log scale apart from the right panel of Figure~\ref{fig:S_processes}.

\begin{figure}
  \begin{center}
    \includegraphics[width=0.48\textwidth]{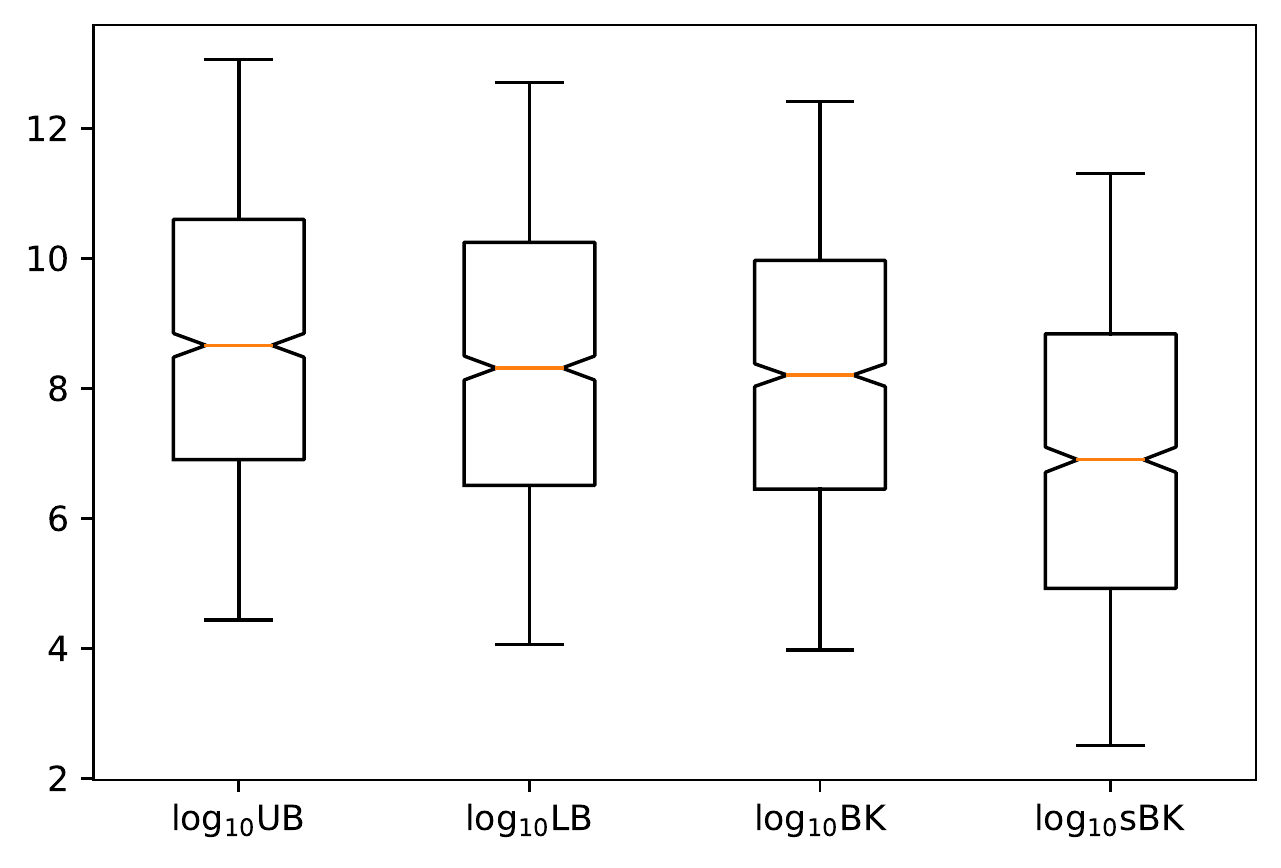}
    \includegraphics[width=0.48\textwidth]{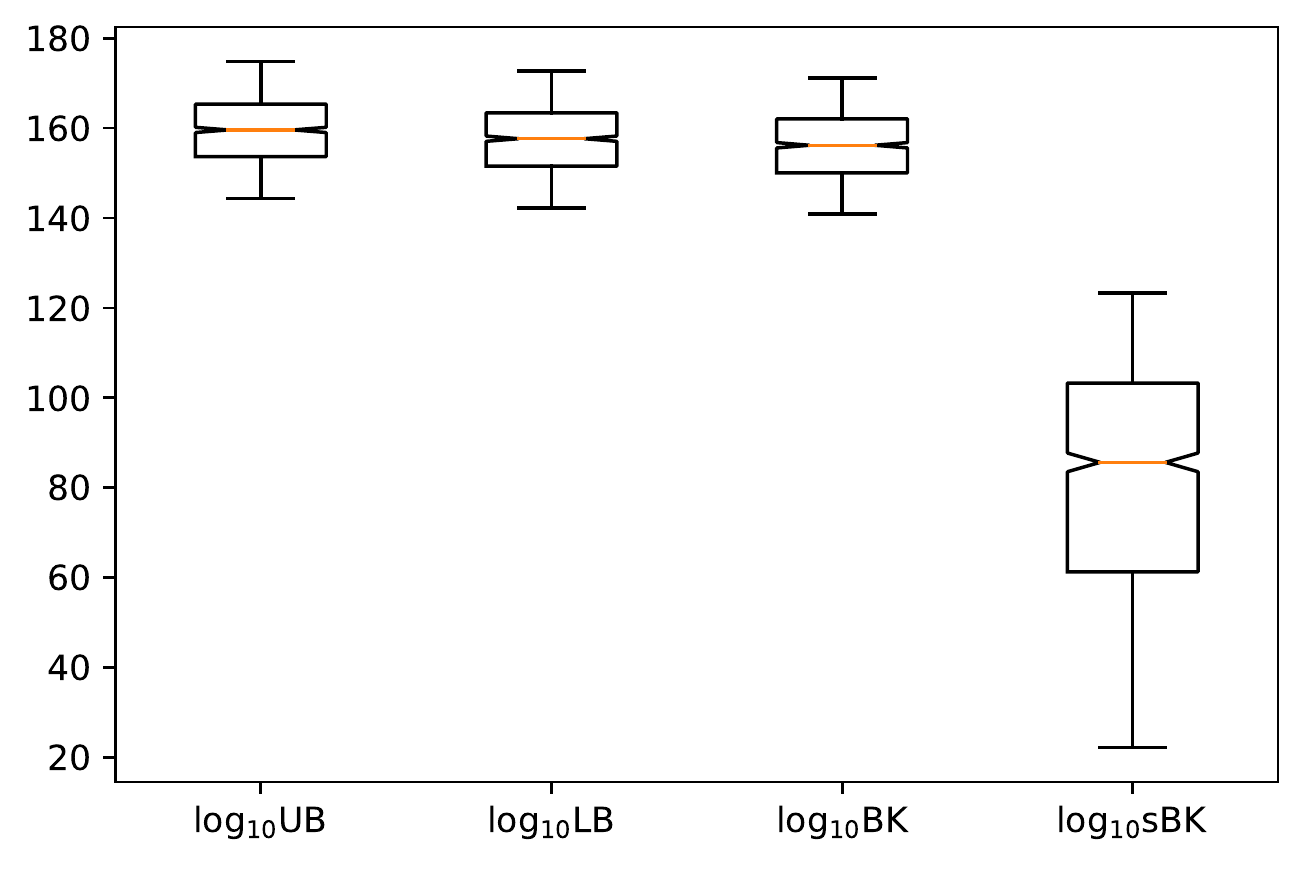}
  \end{center}
  \caption{Boxplots based on $10^3$ runs for the final values of the two benchmarks (upper $\UB$ and lower $\LB$),
    the Bayes--Kelly conformal test martingale (BK), and its simplified version (sBK)
    in the medium scenario.
    Left panel: hard case.
    Right panel: easy case.}
  \label{fig:M_boxplots}
\end{figure}

\begin{figure}
  \begin{center}
    \includegraphics[width=0.48\textwidth]{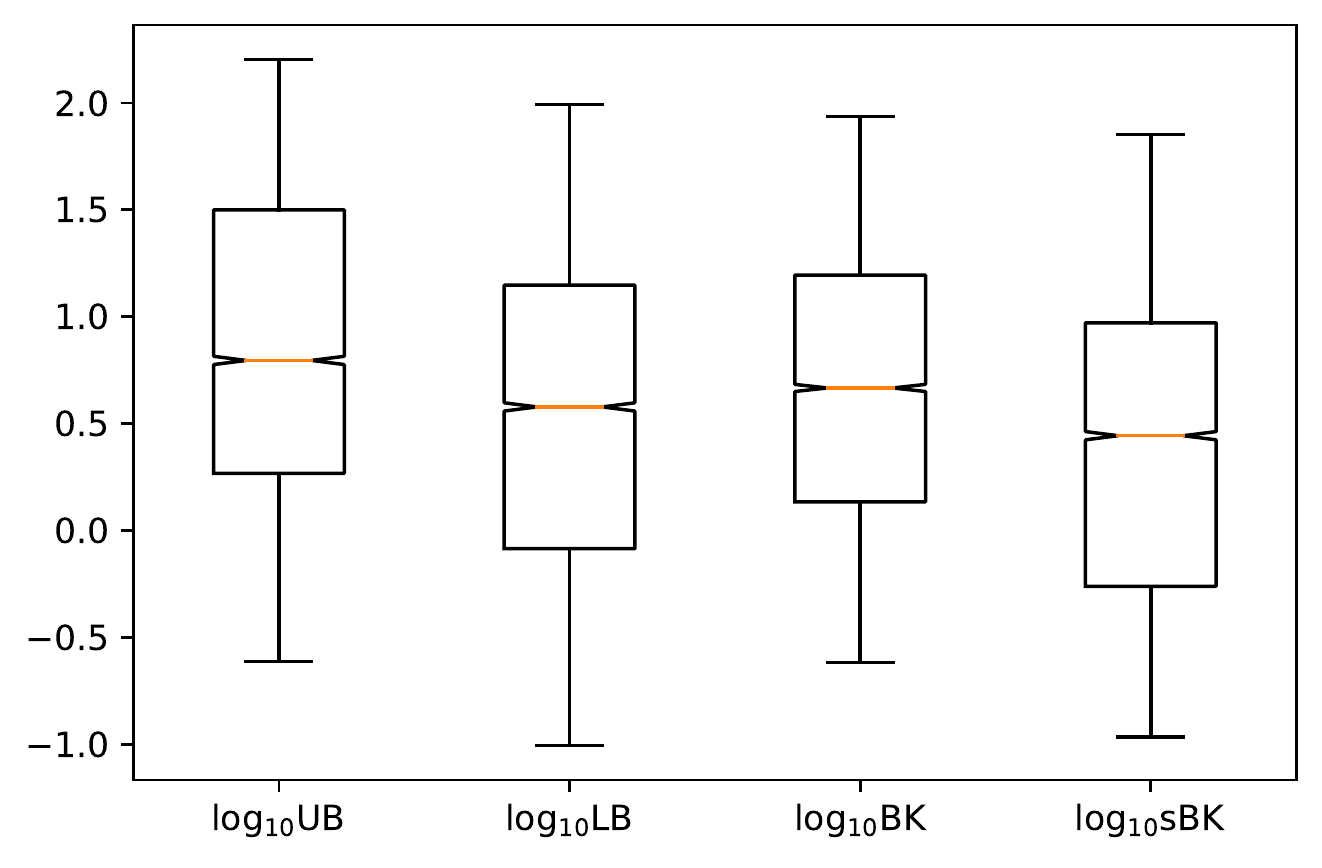}
    \includegraphics[width=0.48\textwidth]{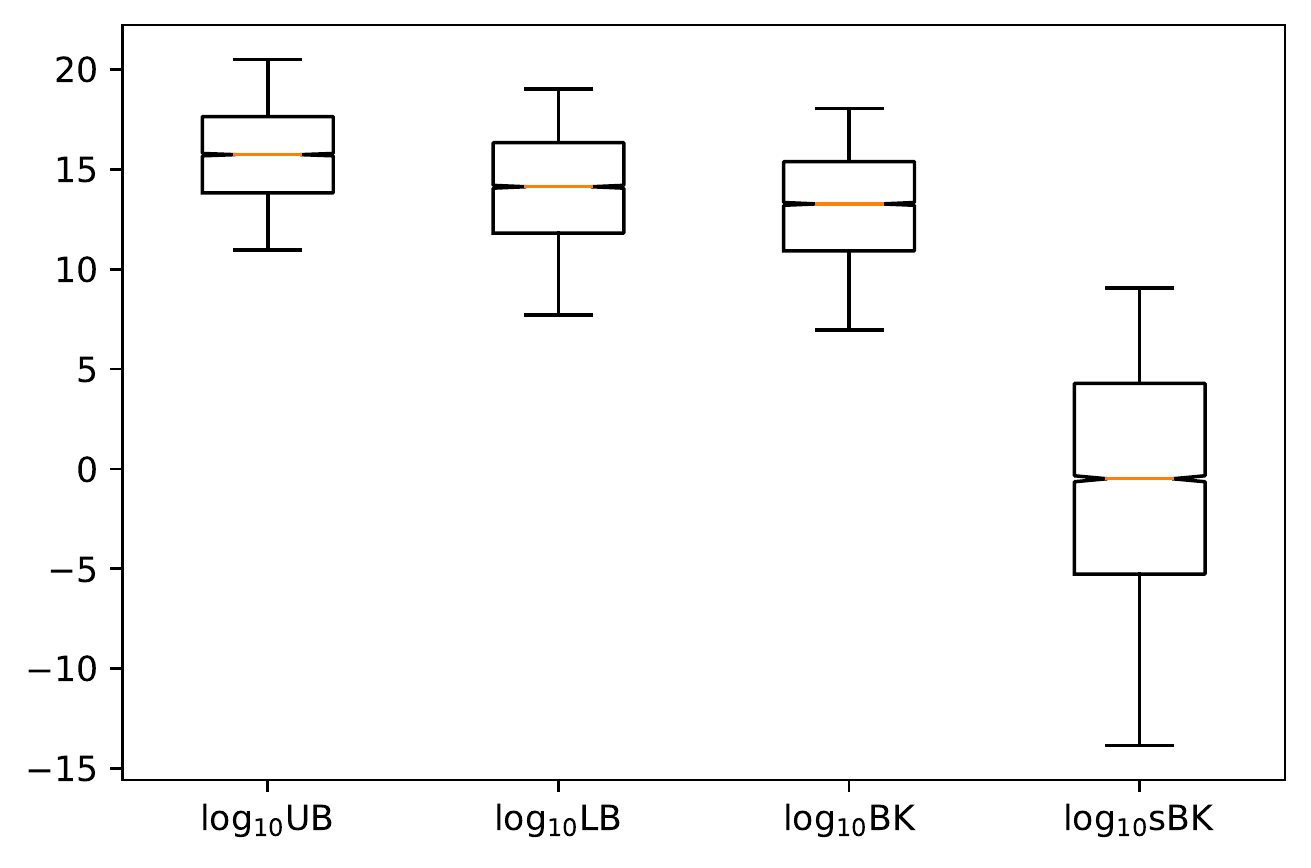}
  \end{center}
  \caption{The analogue of Figure~\ref{fig:M_boxplots} for $10^4$ runs in the small scenario.}
  \label{fig:S_boxplots}
\end{figure}

The right panel of Figure~\ref{fig:S_processes}
and Figures~\ref{fig:M_boxplots} and~\ref{fig:S_boxplots}
show that the statistical performance of the simplified Bayes--Kelly martingale
particularly suffers in the easy case.
The notches in the boxplots in Figures~\ref{fig:M_boxplots} and~\ref{fig:S_boxplots} indicate confidence intervals for the median.

\end{document}